\documentclass[a4paper,12pt,oneside]{article}

\usepackage{graphicx}

\usepackage{latexsym}
\usepackage[english]{babel}
\newtheorem{theorem}{Theorem}
\newtheorem{lemma}[theorem]{Lemma}
\newtheorem{proposition}[theorem]{Proposition}

\newtheorem{definition}[theorem]{Definition}
\newtheorem{corollary}[theorem]{Corollary}
\newtheorem{remark}[theorem]{Remark}

\begin{document}
\title{Information, initial condition sensitivity and dimension in weakly chaotic dynamical systems }

\author{Stefano Galatolo
\\
Dipartimento di Matematica  \\
and \\
 Centro Interdisciplinare per lo Studio dei
Sistemi Complessi,\\
 Universit\`a di Pisa \\
 Via Buonarroti 2/a, Pisa, Italy \\
e-mail: galatolo@dm.unipi.it}

\maketitle
\abstract {We study  generalized  indicators 
of sensitivity to initial conditions and orbit complexity in topological
dynamical systems. The orbit complexity is a measure of the 
asymptotic behavior of the information that is necessary to describe the orbit of a given point. The indicator generalizes, in a certain sense, the Brudno's orbit complexity (which is strongly related to the entropy of the system).
 The initial condition sensitivity indicators consider the  asymptotic behavior of the speed of divergence of nearby staring orbits.
 The indicators have non trivial values also in weakly chaotic 
dynamical systems, characterizing various cases of weakly chaotic
dynamics. Then, using  constructivity, local relations are proved between generalized orbit  complexity,  
initial condition  sensitivity and the dimension of the underlying space or invariant measure. }

\def\stackunder#1#2{\mathrel{\mathop{#2}\limits_{#1}}}

\section{Introduction}

Sensitivity to initial conditions is a feature of chaotic dynamical systems.
A system is sensitive to the initial condition if nearby starting orbits
diverge  during  the evolution.
This is a cause of unpredictability and {\em complexity} of the evolution.

Sensitivity can be measured quantitatively in different  ways. In the
strongly chaotic case,
when a system has exponential sensitivity, i.e. nearby starting 
orbits diverges with exponential speed. 
Two examples are Lyapunov exponents and Brin-Katok local entropy.

As concerns the complexity of the behavior, it has  been measured 
for example by the various notions of entropy and many other invariants
(see for example  \cite{BK}, chapter 3).
In \cite{brud} Brudno  defined the complexity of an orbit of a topological dynamical system, using 
the Kolmogorov (algorithmic) information content of a string (section 2).
 A set of strings is associated by a certain construction to the orbit of a point   and then the complexity of the orbit is defined by the  information content   of the associated strings.
The complexity of an orbit is then a measure of the rate of information 
that is necessary to describe the orbit as the time increases.

 Brudno proved that if the  system is  ergodic  on a  compact space, the entropy of the system  is almost everywhere equal to the orbit complexity. In other words, if such a system  has positive entropy, then for a.e. $x$ the algorithmic information that is necessary to
describe $n$ steps of the orbit of $x$ increases linearly with $n$ and the
proportionality factor is the entropy of the system.
This also implies that if a system has an invariant measure $\mu$, its 
entropy is equal to the mean value of the orbit complexity with respect to $\mu$. Then, in compact dynamical systems Brudno orbit complexity
can be viewed as a pointwise version of entropy.

If  the   sensitivity
is measured by Lyapunov exponents then a relation 
between complexity and sensitivity is the well known
Ruelle-Pesin theorem.  Another relation, which does not 
requires smoothness assumptions, but require the knowledge
of an invariant measure in the definition of the indicator 
of sensitivity is given by the theorem of Brin-Katok \cite{BK}.

The philosophical meaning of the above cited theorems can be resumed as follows.
Under some assumptions the average rate of exponential divergence 
of nearby starting orbits is equal to the average rate of information
 that is necessary to describe the orbits of the points of $(X,T)$.
 In the case of compact dynamical systems  this motivates  
the general Ford's claim \cite{ford} that orbit complexity was a synonym of chaos.
However we remark that if the space is not compact there are examples (\cite{WB}, \cite{io2}) of systems
with no sensitivity to initial conditions  and hight orbit complexity.
 
One of the consequences of the results proved in section 6 is that in the compact case the assumption of constructivity recovers a link between complexity 
and sensitivity, both in the positive, and in the zero entropy case.

In many  examples of dynamical systems  the entropy is 0, the speed of
separation of nearby starting trajectories is less than
exponential and the increasing of the information contained in $n$ step of the orbit could be less than linear. This is the case of the so called {\em Weakly Chaotic Dynamics}.

The question naturally arises if a quantitative relation between 
complexity and sensitivity can be found for zero entropy
dynamical systems. This is the main aim of the paper \cite{ionew}.
Our answer is that if the system is constructive, that is, if 
the transition map can be approximated at a given accuracy (in a sense
that will be clarified in section 2) by some algorithm then we have 
such a quantitative relation.
The relation implies for example that stretched exponential initial
condition sensitivity \footnote{If two points starts at distance $\Delta x(0)$ then $\Delta x(t)\sim\Delta x(0) 2^{Ct^{\alpha}}$.} implies power law 
behavior of the quantity of information that is necessary to describe
an orbit.
In this paper we improve the results of \cite{ionew}  defining the indicators
by simpler
notations. This gives to the statements a simpler and stronger form, involving connections with the concept of dimension.


The study of weakly chaotic dynamics was discovered to be
important for application purposes.
There are connections with many physical and economic phenomena: self organized criticality, the so called chaos
threshold, the anomalous diffusion processes  and many others.
In these examples of weakly chaotic dynamics the traditional indicators
of chaos (K.S. entropy, Lyapunov exponents, Brudno's orbit complexity)
vanishes. These indicators are not able to distinguish between all the various
cases of weakly chaotic dynamics. 

Some definitions of generalized entropy have been already proposed in 
literature (see e.g. \cite{takens},\cite{meson},\cite{tsallis} for definitions about dynamical systems). The most fruitful ones were given by  Reniy and Tsallis  that found a great
variety of application. See for example
http://tsallis.cat.cbpf.br/biblio.htm for an updated bibliography of topics related to Tsallis entropy.
In the literature relations between Tsallis entropy and initial data sensitivity have been  proved (\cite{tsallis},\cite{jin}),
indicating that the main field of application of Tsallis entropy in
dynamical systems is the case of power law initial data sensitivity (if two points starts at distance $\Delta x(0)$ then $\Delta x(t)\sim\Delta x(0) t^{\alpha}$).

In what follows it comes out that the two ingredients we need to add to have interesting local statements are constructivity and dimension. 
To give a formal approach to constructivity and algorithms ``acting'' on metric spaces we introduce the notion of computable structure.
The formal definition and some heuristic motivations of computable structure
will be given in section 2, after a short introduction to the notion of algorithmic information content. 

In section 3 we introduce the generalized indicators of orbit complexity.
These indicators characterizes the asymptotic behavior as $n$ increases of the quantity of information that is necessary to describe $n$ steps of the orbit of a point.
We give two definitions of such an indicator.
The first is a straightforward generalization of the Brudno's definition of
orbit complexity. The second is a generalization of the definition
that was given in \cite{io2} to extend Brudno's definition to the non compact case.
In the compact case the two definitions are equivalent.
The second definition is related to the notion of computable structure and 
allows to use the features of constructivity. 

In section 4 we give definitions
of initial conditions sensitivity indicators.
The definition is similar to the one given by Brin and Katok \cite{BK} for the local entropy, but does not use the invariant measure. and depends only on the topological (and metric, in the sense of distance) features of the dynamics.

In section 5 we introduce a notion of complexity of points in a metric space with a computable structure. Such a complexity is a measure of the quantity
of information that is necessary to approximate a point of the metric space at a given accuracy. This notion is related with the notion of dimension, as it was already remarked in \cite{io1}.

In section 6 we state the relations linking orbit complexity, sensitivity and 
dimension. Under the assumption of constructivity we prove local statements
of the form
$$ dimension \times indicator \ of \ min \ sensitivity \ at \ x \leq orbit \ complexity \ at \ x $$

$$ orbit \ complexity \ at \ x \leq dimension \times  indicator \ of \ max \ sensitivity \ at \ x $$

\noindent we remark that these statements are obtained only under the assumption of constructivity (and continuity) of the transition map. No differentiability or hyperbolicity
assumptions are needed. The equations gives information both in the positive and null entropy case.


\section{Algorithmic Information Theory and Computable Structures}

%

In this section we give a short  introduction to algorithmic information theory and to computable structure and constructivity.
A more detailed exposition of algorithmic information theory can be found in \cite{Zv} or \cite{Ch}. 

Let us consider the set  $\Sigma=\{0,1\}^*$ of finite (possibly empty)  binary strings. If $s$ is a string we define $|s|$ as the length of $s$.

Let us consider a Turing machine (a computer) $C$: by writing $C(p)=s$ we mean that $C$ starting with input $p$ (the program) stops with output $s$ ($C$ defines a  partial recursive function  $C:\Sigma \rightarrow \Sigma$). If the input gives a never ending computation the output (the value of the recursive function) is not defined.
If $C:\Sigma \rightarrow \Sigma$ is recursive  and its value is defined
for all the input strings in $\Sigma$ (the computation stops for each input)
then we say that $C$ is a {\em total} recursive function from $\Sigma$ to $\Sigma$. The algorithmic information content (A.I.C.) of a string will be the length of the shortest
program that outputs the string.

\begin{definition}
The Kolmogorov complexity or Algorithmic Information Content of a string 
$s$ given $C$ is the length of the smallest
program $p$ giving $s$ as the output:

$$
{ AIC}_C(s)=\stackunder{C(p)=s}{\min }|p| 
$$

\noindent  if $s$ is not a possible output for
the computer $C$ then ${AIC}_C(s)=\infty $ .
\end{definition}
Let us see some examples.
The algorithmic information content 
of a $2n$ bits   long {\em periodic} string 
$$s=''1010101010101010101010...''$$  is small because 
the string is output of a shortest program: 

\centerline{\em repeat $n$ times (write (``10''))}

\noindent the AIC of the string $s$ is then less or equal than $ log(n)+Constant$ where $log (n)$ bits are sufficient  to  code ''$n$'' and the constant represents the length of the code for the computer $C$ representing the instructions ''repeat, write...''. As it is intuitive the information content of a periodic string is very poor.
On the other hand each, even random  $n$ bits long string $$s'=''1010110101010010110...''$$ is output of the trivial program 
$${write(''1010110101010010110...'')}$$ 
this is of length $n+constant$. This implies that the AIC of each string
is (modulo a constant which depends on the chosen computer $C$) less or equal than its length.

From the last definition  the algorithmic information content of a string depends on the choice of the computer $C$. There is a class of computers that 
allows an ''almost'' universal definition of algorithmic information content of a string: if we consider computers from this class the A.I.C. of a string
is  defined independently of the computer up to a constant. In order
to define such a class of universal computers 
 we give some notations that are necessary to work with strings:
let us consider   $c:\Sigma \rightarrow {\bf N}$ from the set $\Sigma$  and the set \-${\bf N}$ associating strings and natural numbers in the following 
natural way 

$$ \emptyset \rightarrow 0, \ \ \  0 \rightarrow 1 , \ \ \  1
\rightarrow 2, \ \ \  00 \rightarrow
3, \ \ \  01\rightarrow 4 , \ \ \ 10\rightarrow 5 , ...$$

 This correspondence  allows us to interpret  natural numbers as  strings
and vice-versa when it is needed. 
 We remark  that
$|s|\leq log(c(s))+1$. \footnote{ In this paper all the logarithms are in base two.}
If  $s$ is a string with $|s|=n$ let us denote by  $s\hat{\ \ } $ the string $s_0s_0s_1s_1...s_{n-1}s_{n-1}01.$ If $a=a_1...a_n$ and $b=b_1... b_m$ are strings
then $ab$ is defined as the string $a_1...a_nb_1... b_m$. If $a$ and $b$ are strings
then $a\hat{\ \ } b$ is an encoding of the couple $(a,b)$. There is
an algorithm that getting the string $a\hat{\ \ } b$ is able to recover both the strings $a$ and $b$. 
An universal Turing machine intuitively is a machine that can emulate 
any other  Turing machine if an appropriate input is given. Before to give
the formal definition 
we recall that there is a recursive enumeration $A_1,A_2...$
of all the Turing machines.    
\begin{definition}
A  Turing machine ${\cal U}$ is said to be universal if for all $m\in {\bf N}$ and $p\in \Sigma$ then  ${\cal U}(c^{-1}(m)\hat{\ \  } p)=A_m(p)$.

\end{definition}

In the last definition the machine ${\cal U}$ is universal because ${\cal U}$ is able to emulate each other machine
$A_m$ when in its input we specify the number $m$ identificating $A_m$ and
the program to be runned by $A_m$.
It can be proved  that an universal
Turing machine 
exists.

\begin{definition}
A  Turing machine $ F$ is said to be {\em asymptotically optimal }if for each Turing machine $C$ and each binary string $s$ we have $AIC_{F}(s)
\leq {AIC}_C(s)+c $ where
the constant $c$ depends on $C$ and not on $s$.
\end{definition}

The following proposition can be proved from the definitions

\begin{proposition}
If ${\cal U}$ is an universal Turing machine
then ${\cal U}$ is asimptotically optimal.
\end{proposition}
 
 This tells us that choosing an universal Turing machine the
complexity of a string is defined independently of the given 
Turing machine  up to a constant.  
 For the remaining part of the
paper we will suppose that  an universal Turing machine ${\cal U} $ is chosen once forever.

\subsection{Computable  Structures, Constructivity }

Many models of the real world use the notion of real numbers or more
in general the notion of complete metric spaces. Even if you consider a very
simple complete metric space, as, for example, the interval $\left[
0,1\right]  $ it contains a continuum of elements. This fact implies that most
of these elements (numbers) cannot be described by any finite alphabet.
Nevertheless, in general, the mathematics of complete metric spaces is simpler
than the ''discrete mathematics'' in making models and the relative theorems.
On the other hand the discrete mathematics allows to make computer
simulations.
This is one of the reasons why we  introduce the
notion of {\em computable structure} which is a way to relate the world of
continuous models with the world of computer simulations. 
The set of binary strings $\Sigma$ is the mathematical
abstraction of the world of the ''computer'', or more in general is the
mathematical abstraction of the ''things'' which can be expressed by any
language. The real objects which we want to talk about are
modeled by the elements of a metric space $(X,d).$ We will  interpret the
objects of $\Sigma$ as points of $X.$  A
 computable structure on a separable metric space $(X,d)$ is a class of dense {\em interpretations } ($I:\Sigma\rightarrow X$)  of the
space of finite strings $\Sigma$ in the metric space. The {\em interpretations}  are such that the distance
$d$ restricted to the points that are images of strings ($x=I(s):x\in X,s\in\Sigma$) is a ''computable'' function.
The use of computable structures allows to consider algorithms ''acting''
over metric spaces and to define constructive functions between metric spaces,
that is, functions such that we can work with by using a finite amount 
of information. Many concrete metric spaces used in analysis or in geometry have a
natural choice of a computable structure.
 In the following we  will assume that the
dynamical systems under our consideration are constructive. All the dynamical systems that we can construct explicitely 
are construcive.
From the philosophical point of wiew we  think  that the assumption of constructivity is not unnatural because
 even if the  maps coming from  physical reality were not constructive, the models used  to describe such a reality should be constructive (to allow calculations).
On the other hand, to add constructivity  allows to prove stronger
theorems, avoiding pathologies coming from random maps.

As it was said before
an interpretation function is a way to interpret  a string as a point of the metric space.

\begin{definition}
An interpretation function on  $(X,d)$ is a function $I:\Sigma \rightarrow
X$ such that $I(\Sigma)$ is dense in $X$.
\end{definition}

 A point $x\in X$ is said to be {\em ideal} if it is the image of some string $x=I(s), s\in \Sigma$.  
 An interpretation is said to be computable  if the distance between ideal points is computable with arbitrary precision:

\begin{definition}
A computable  interpretation function on  $(X,d)$ is a function $I:\Sigma \rightarrow
X$ such that $I(\Sigma)$ is dense in $X$ and  there exists a total recursive
function $D:\Sigma \times \Sigma \times
 {\bf N} \rightarrow {\bf Q}$ such that $\forall s_1,s_2 \in \Sigma ,n\in
{\bf N}$:

$$|d( I(s_1),I(s_2))-D(s_1,s_2,n)|\leq \frac{1}{2^n}.$$

\end{definition}

Two interpretations are said to be equivalent if the distance from an ideal
point from the first and a point from the second is computable up to
arbitrary precision.
For example, the finite binary  strings $s\in \Sigma $ can be interpreted as rational numbers by interpreting the string as the binary expansion of a number. Another interpretation can be given by
interpreting a string as an encoding of a couple of integers whose ratio
gives the rational number. If the  encoding is recursive, the two interpretation are equivalent.

\begin{definition}\label{def7}     
Let $I_1$ and $I_2$ be two computable interpretations in $(X,d)$;
we say that $I_1$ and $I_2$ are equivalent if there exists a total recursive function
$D^*:\Sigma \times \Sigma \times {\bf N} \rightarrow {\bf Q}$,
such that $\forall s_1,s_2 \in \Sigma ,n{\in {\bf N}}$:

$$|d( I_1(s_1),I_2(s_2))-D^*(s_1,s_2,n)|\leq \frac{1}{2^n}.$$

\end{definition}

\begin{proposition}
The relation  defined by definition \ref{def7} is an equivalence relation.
\end{proposition}

\noindent For the proof of this proposition see \cite{io2}.

\begin{definition}
A computable structure ${\cal I}$ on $X$ is an equivalence class  of computable  interpretations  in $X$.
\end{definition}

For example if $X={\bf R}$ we can consider the interpretation $I:\Sigma \rightarrow {\bf R}$ defined in the following way: if $s=s_1...s_n\in \Sigma$ then \begin{equation}\label{I}I(s)=\sum_{1\leq i \leq n}s_i2^{[n/2]-i}.\end{equation} This is an interpretation of a string
as a binary expansion of a number.
  $I$ is a computable interpretation, the computable structure on ${\bf R}$ containing $I$ will be called {\em standard} computable structure.
If $r=r_1r_2...$ is an infinite string such that $lim \frac{K_C(r_1...r_n)}{n}=1$\footnote{such a string exist, see for example \cite{io1} theorem 13.} then the interpretation $I_r$ defined as $I_r(s)=I(s)+\sum r(i)2^{-i}$ is computable but not equivalent to $I$. $I$ and $I_r$ belongs to different computable structures.
 
In a similar way it is easy to construct computable structures in ${\bf R}^n$ or in separable function spaces codifying a dense subset (for example the set of step functions) with finite strings. 
We remark as a property of the computable structures that if $B_r(I(s))$ is 
an open ball with center in an ideal  point $I(s)$ and rational radius $r$ and $I(t)$ is another point   then there is an algorithm that 
verifies if  $I(t)\in B_r(I(s)) $. If  $I(t)\in B_r(I(s)) $ then the algorithm
outputs ''yes'', if $I(t)\notin B_r(I(s)) $ the algorithm  outputs ''no'' or 
does not stop. The algorithm calculates $ D(s,t,n)$ for each $n$ until it
finds that $D(s,t,n)+2^{-n}< r$ or  $D(s,t,n)-2^{-n}> r$, in the first case
it outputs ''yes'' and in the second it outputs ''no'', if $d(I(s),I(t))\neq r$ the algorithm will stop and output an answer.

We give a definition of {\em morphism} of metric spaces with computable 
structures, a morphism is heuristically  a computable
function between  computable metric spaces.

\begin{definition}
If $(X,d,{\cal I})$ and $(Y,d',{\cal J}) $ are spaces with  computable
structures; a function $\Psi :X\rightarrow Y$ is said to be  a morphism of
computable structures if $\Psi $ is uniformly continuous and for each pair
$I\in {\cal I},J\in {\cal J}$ 
there exists a total recursive function
$D^*:\Sigma \times \Sigma \times {\bf N} \rightarrow {\bf Q}$,
such that $\forall s_1,s_2 \in \Sigma ,n{\in {\bf N}}$:

$$|d'(\Psi( I(s_1)),J(s_2))-D^*(s_1,s_2,n)|\leq \frac{1}{2^n}.$$
\end{definition}

\noindent We remark that $\Psi $ is not required to have dense image
and then $\Psi(I(*))$ is not necessarily an interpretation function 
equivalent to $J$.

\begin{remark}\label{remark10} 
As an example of the  properties of the morphisms, we remark that if a map $\Psi:X\rightarrow Y $  is a morphism  then given a point $x\in I(\Sigma) \subset X$ it is possible to find by an algorithm a point $y\in J(\Sigma) \subset Y$ as near as we want to $\Psi (x) $. 
\end{remark}
The procedure is simple: if $x=I(s)$ and we want to find a point $y=J(z_0)$ such that $d'(\Psi( I(s)),y)\leq 2^{-m} $ then we calculate $D^*(s,z,m+2)$ for each $z\in \Sigma$ until we find $z_0$ such that $D^*(s,z_0,m+2)<2^{-m-1}$. Clearly $y=J(z_0)$ is such that $d'(\Psi(x),y)\leq 2^{-m} $. The existence of such a $z_0$ is assured by the density of $J$ in $Y$.
In particular the identity is a morphism.

We also remark that by a similar procedure, given a point $I(s_0)$ and $\epsilon \in {\bf Q}$ it is possible to find a point $I(s_1)$ such that $d(I(s_0),I(s_1))\geq \epsilon $.


A constructive map is a morphism for which the continuity relation between $\epsilon$ and $\delta$  is given by a
 recursive function. The following is in some sense a generalization of the definition of Grzegorczyk, Lacombe (see e.g. \cite{purel}) of constructive function.

\begin{definition}\label{definition11}
A function $\Psi:(X,d,{\cal I})\rightarrow (Y,d',{\cal J}) $ between  spaces with computable structure
$ (X,{\cal I}) ,(Y,{\cal J})$ is said to be constructive
if  $\Psi$ is a morphism between the  computable structures and it is 
effectively uniformly continuous, i.e. there is a total recursive function
$f:{\bf N}\rightarrow {\bf N} $ such that for all $x,y \in X$  
$d(x,y)<2^{-f(n)}$ implies $d'(\Psi (x),\Psi (y))<2^{-n}$.
\end{definition}

If $X$ is a space with a computable structure and $T:X\rightarrow X$ is constructive then we call $(X,T)$ a constructive dynamical system.

The following Lemma states that
if a map between spaces with a computable structure is constructive then there
is an algorithm to follow the orbit  each ideal point $x=I(s_0)$.
 The proof can be found in \cite{ionew}

\begin{lemma}\label{lemma12}
If $T: (X,{\cal I})\rightarrow (X,{\cal I})$  is constructive, $ I\in{\cal I}$ then  there is an algorithm (a total recursive function)
$A:\Sigma \times {\bf N}\times {\bf N}\rightarrow \Sigma$ such that $\forall k,m\in {\bf N},s_0\in \Sigma$ $d(T^k(I(s_0)),I(A(s_0,k,m)))< 2^{-m} $. 
\end{lemma}

\section{Brudno's definition of orbit complexity and its generalization}

Here we sketch Brudno's definition of orbit complexity.
The construction will be generalized, to provide
invariants that are useful in the weakly chaotic case.

 Let us consider a topological  dynamical system $(X,T)$.  $X$ is a 
metric space and $T$ is a continuous \footnote{However we remark that for the next definition to be well posed the continuity is not strictly necessary.}
  mapping  $X\rightarrow X$. Let us consider a  finite open cover  $\beta=\{B_0,B_1,...,B_{N-1}\}$ of $X$. We use this cover
to code the orbits of $(X,T)$  into a set of infinite strings. 
A symbolic coding of the orbits of $X$ with respect to the open cover $\{B_i\} $ is a string listing the sets $B_1,..,B_n$ visited by the orbit of $x$ during 
the iterations of $T$. Since the sets $B_i$ may have non empty intersection 
then an orbit can have more than one possible coding. More precisely.
If  $x\in X$ let us define the set of symbolic orbits of $x$ with respect to $\beta$ as: 
$$\varphi _\beta(x)=\{\omega \in \{0,1,...,N-1\}^{\bf N}:\forall n\in
{\bf N} , T^n(x)\in B_{\omega (n)}\}.$$

\noindent The set $\varphi _\beta(x)$ is the set of all the possible 
codings of the orbit of $x$ relative to the cover $\beta$. 

The Brudno's definition of orbit complexity of $x$ with respect to $\beta $
is $$K(x,T,\beta)=\stackunder {n\rightarrow \infty}{limsup}
\stackunder {\omega \in \varphi _\beta (x)}{min}\frac {AIC_U(\omega ^n)}{n}
$$
\noindent where $\omega ^n$ is the string containing the first $n$ digits of $\omega $. We remark that $\omega  ^n$ is not a binary string. It is easy to imagine how the definition of algorithmic information content  can be extended to strings made of digits coming from a finite alphabet. 
 This definition measures the average quantity of information that is necessary to describe a step of the orbit of $x$ by the open sets of $\beta$. The Brudno's  orbit complexity is positive when the quantity of information increases linearly. This is the case of positive entropy (see theorem \ref{15}).
When the system is weakly chaotic the information increases less
than linearly and  Brudno's complexity is 0. As we will see next there
are many possible different behavior of the quantity of information
in weakly chaotic dynamical systems. The generalized indicator
will be able to distinguish between these different behaviors.

We  give a measure of such an  asymptotic behavior by comparing the quantity of information necessary to describe $n$ step of the orbit  with a  function $f$  whose asymptotic behavior is known. For each monotonic function $f(n)\rightarrow \infty$ we define an indicator
of orbit complexity by comparing the asymptotic behavior of  ${AIC}_U(\omega ^n)$ with $f$.
The complexity of the orbit of $x\in X$ relative to $f$ and  $\beta$ is defined as:

$$K^f(x,T,\beta)=\stackunder {n\rightarrow \infty}{limsup}
\stackunder {\omega \in \varphi _\beta (x)}{min}\frac {AIC_U(\omega ^n)}{f(n)}
$$

 Taking the supremum   over the set of all {\em finite open} covers
$\beta$ of the metric space $X$
 it is possible to define the complexity of the orbit of
$x$: 

$$K^f(x,T)=\stackunder {\beta}{sup}( K^f(x,T|\beta))$$

This definition associates to a point belonging to $X$ and a function $f$ a real number which is a measure of the complexity of the orbit of $x$ with respect to the asymptotic behavior of $f$. For example, if $f$ is the identity: $f(n)=n$ this definition coincides with the original Brudno's one.
We remark that it is important to suppose that sets in the covers are open
(see the discussion in \cite{brud}, \cite{WB},\cite{io2},\cite{Licata}).

Generalized orbit complexity is invariant under topological conjugation,
as it is stated in the following theorem whose proof follows directly from the definitions:

\begin{theorem} If the dynamical systems $(X,T)$ and $(Y,S)$ are topologically  conjugate, and  $\pi:X\rightarrow Y$ is the conjugating homeomorphism, and
$\pi(x)=y$  then $K^f(x,T)=K^f(y,S)$.
\end{theorem}

In the literature (see e.g. \cite{io2}, \cite{brud}, \cite{Wh}, \cite{WB}, \cite{Gaspard},\cite{keller})   many relations have been proved between orbit complexity and
other forms of complexity of a dynamical system (Kolmogorov entropy, topological entropy and others) and with 
other problems concerning  orbits of a dynamical system. 
The following is of particular interest:

\begin{theorem}\label{15}(Brudno's main theorem.)
Let $(X,T)$ be a dynamical system over a compact space. If $\mu$ is an ergodic probability
measure on $(X,T)$, then

$$K^{id}(x,T)=h_{\mu}(T)$$
 
\noindent for $\mu$-almost each $x\in X$.
\end{theorem}

Where $h_{\mu}(T)$ is the Kolmogorov entropy of $(X,T)$  with respect
to the invariant measure $\mu$.
We also recall the following lemma which will be used lather.
\begin{lemma}\label{16}
If $\alpha $ and $\beta $ are open covers of $X$ and $\alpha $ is a refinement of $\beta$ then 
$$K^f(x,T,\beta)\leq K^f(x,T,\alpha)$$

\end{lemma}

From now on in the notation $K^f(x,T)$ we will avoid to explicitly mention the map $T$ when it is clear from the context. 
We now give some example of different behaviors of $K^{f}(x)$.
If $x$ is a periodic point it is easy to see that since the information content of a $n$ digit long periodic string is $\leq log(n)+C$ then  ${K}^{f}(x)=0$ if $log(n)=o(f)$, moreover ${K}^{log(n)}(x)=1$, this follows by the (not so trivial)  remark
that if $\omega$ is a  periodic string then $\stackunder {n\rightarrow \infty}{ limsup} \frac {AIC(\omega ^n)}{log (n)}=1$ (while the $lim inf$ is  equal zero). 
By theorem \ref{15} it follows  that if a system is compact, ergodic and has positive Kolmogorov entropy then for almost all points we have ${K}^{f}(x)=\infty$ if $f=o(id)$ and ${K}^{n}(x)=h_\mu$.
We also remark that (when the space is compact)  the linear one is the maximum over all the possible asymptotic behaviors.
Indeed if $X$ is compact, for each $\epsilon$ there is a finite cover $V$ made
of balls with radius $\epsilon$, then a program that follows $n$ steps of
the orbit of any point with the accuracy $\epsilon$ can be simply made by  listing $n$
balls of the cover, then, if $X$ is compact then $
\stackunder {\omega \in \varphi _\beta (x)}{min}{AIC_U(\omega ^n)}
\leq \alpha n+C$, where $\alpha $ is the logarithm of the number of elements in the cover $V$.
Another important example is the piecewise linear Manneville map:
\begin{equation}\label{**}
T_z(x)=\left\{ \begin{array}{cc}
\frac{\xi _{k-2}-\xi _{k-1}}{\xi _{k-1}-\xi _k}(x-\xi _k)+\xi _{k-1} & \xi
_k\leq x<\xi _{k-1} \\ 
\frac{x-a}{1-a} & a\leq x\leq 1
\end{array}
\right. 
\end{equation}

\noindent with $\xi _k=\frac a{(k+1)^{\frac 1{z-1}}}$ ,$\ k\in N,z \in {\bf R}, \ z\geq 2$.
This is a P.L. version of the Manneville map $T(x)=x+x^z \ (mod\  1)$ ({\em see fig 1}).

\begin{figure}
\begin{center}
\includegraphics[height=7cm]{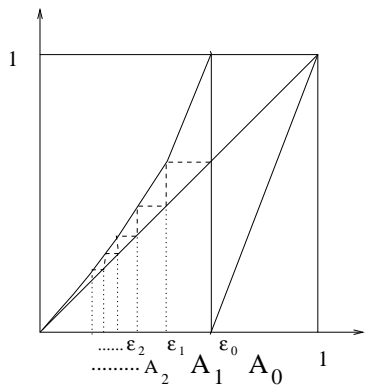}
\caption{\label{xsosx}The  map  $T_z$ and the partition $A_i$}
\end{center}
\end{figure}
 The Manneville map was introduced in \cite{Manneville}  as an extremely simplified 
model of intermittent behavior in turbulence, then its mathematical 
properties was studied by many authors (e.g.  \cite{Gaspard},\cite{pollicot}). And the map was applied as a model of other physical phenomena ( \cite{Grigolini},\cite{Grigolini96}).

It is proved in \cite{ionew} and \cite{bonanno} following the ideas of \cite{Gaspard}  that for almost each $x$ (for the Lesbegue measure) 
${K}^{f}(x)=0$ if $n^\alpha=o(f(n)) $, with $\alpha=\frac 1{z-1}$ . And  ${K}^{f}(x)=\infty$ if $f(n)=o(n^\alpha) $. Then for almost each point we have that the information content
of an orbit is Lesbegue a.e. increasing as a power law with exponent $\alpha$.

\subsection{Another definition of  orbit complexity}
We now give another definition of orbit complexity.
This definition is equivalent to the previous one  when the space 
is compact and has a computable structure. The new definition
 allows 
to exploit the features of constructivity.
Another feature of the new definition is that extends Brudno's definition
to the non compact case (see, e.g. \cite{io2}).

Let $(X,{\cal I})$ a separable metric space with a computable structure
${\cal I}$.
 To interpret the output of a calculation which is a finite string as a finite sequence in $X$ let us consider an interpretation function $I$ and 
 a total recursive surjective function.  
$${\cal Q}:\Sigma \rightarrow \Sigma ^* $$
where $\Sigma ^*$ is the set of finite sequences in $\Sigma $.
${\cal Q}$ associates to a string a sequence of strings that will be 
interpreted as a sequence of points as follows.
 Now let us
consider an  universal Turing machine ${\cal U}$.  For each program $p$  we  define $U(p)\in X^* $ (the set of finite sequences in $X$) as 

$$U(p)={I}( {\cal Q}({\cal U}(p)))$$

\noindent where $I$ is extended in the  obvious way  to a map from the space $\Sigma ^*$ to $X^*$. $U_i(p)\in X$ is defined as the $i-$th point of $U(p)$ .
With this definition we can interpret the output of a calculation
as a finite sequence in $X$.
We  remark that given $\cal Q $ and a sequence of strings $s_1,...,s_n$
it is possible by an algorithm to    find a single string $s$ such that ${\cal Q} (s)=(s_1,...,s_n)$.

\begin{definition}
We define the algorithmic information content of the sequence $x,T(x),...,T^{n}(x)\in X^*$ at accuracy $\epsilon $ and with respect to the interpretation $I$
as:

$${\cal K}^I(x,n,\epsilon) =\min \left\{ |p|\ s.t. \ U(p)\in X^{n+1},
\stackunder {0\leq i\leq n} {max} (d(U_{i}(p),T^i(x)))<\epsilon\right\}. $$  
\end{definition}

\noindent The calligraphic ${\cal K} $ is used to distinguish this 
definition from the definition of Brudno.

 As before we  choose a function $f$ and consider the asymptotic behavior for $n\rightarrow \infty$ and define ${\cal K}^{f,I}(x,\epsilon) :X\times {\bf R}\rightarrow {\bf R}$ as:
$${\cal K}^{f,I}(x,\epsilon)=\stackunder {n\rightarrow \infty}{limsup}\frac {{\cal K}^I(x,n,\epsilon)}{f(n)}.$$

\begin{remark}$ {\cal K}^{f,I}(x,\epsilon)$ is a non increasing function
with respect $\epsilon$.
\end{remark}

Finally, we consider
the behavior when $\epsilon$ goes to $0$ and   we define $ { \cal K}^{f,I}(x) :X\rightarrow {\bf R}$ as

\begin{definition}
The orbit complexity of $x$ with respect to $f$ and the interpretation $I$ is defined as:
$${\cal K}^{f,I}(x)=\stackunder {\epsilon\in {\bf R}^+}{sup}   {\cal K}^{f,I}(x,\epsilon) .  $$
\end{definition}

\begin{theorem}\label{teorema2}
 If $X$ is compact  then for each computable interpretation ${ I}$, for each $x\in X$  
$${\cal K}^{f,{I}}(x)=K^f(x)$$
\end{theorem} 

This result shows that in the compact case the 'calligraphic' orbit
complexity does not depends  on the choice of the computable
structure on the space.

{\em Proof.} 
{\em First part.}  Let us consider an open cover $V$ and $K^f(x,V)$, if $V_\epsilon$ is a refinement  of $V$ made of balls $B_0,...,B_m$ with radius $\epsilon $, by Lemma \ref{16} we have 
$$K^f(x,V_\epsilon)\geq K^f(x,V)$$ 
 since each cover has a refinement made of $\epsilon$ balls,
then $$K(x)=\stackunder {V_\epsilon \in \epsilon -covers} {sup} K(x,V_\epsilon).$$ Let us consider $I $ and consider  $s_0,...,s_m$ such that $B(I{(s_j)},2\epsilon)\supset B_j$
 let us prove that for each $V_\epsilon $
\begin{equation}\label{star} K^f(x,V_\epsilon)\geq {\cal K}^{f,I}(x,2\epsilon)
\end{equation}
this follows from the fact that there is a constant $C$ such that 
$K(x,n,V_\epsilon)+C\geq {\cal K}^{I}(x,n,2\epsilon)$.
Indeed given  a minimal program $p$ generating $n+1$ steps of a symbolic orbit 
of $x$ with respect to $V_\epsilon$: ${\cal U}(p)=\omega^n=(i_0,..,i_n)$
there is a  a program $p'$ containing $p$ that codifies the following procedure

1) it calculates $\omega ^n$ running the given program $p$ 

2) it has a finite list that associates to each symbol in $\omega ^n$
representing a ball $B_i$ of the cover $V_\epsilon $ the
string $s_i$.
With this list it associates to $\omega ^n$ a finite sequence of strings $
s_{i}$

3)it calculates the single string $s^{*}$ such that ${\cal Q}(
s^{*})=(s_1,...,s_n), $ and output this string. Then $ {\cal K}^f(x,n,2\epsilon)\leq |p'|, |p'|=|p|+C$, then $K^f(x,V_\epsilon)\geq {\cal K}^{f,I}(x,2\epsilon)$  and then the statement is proved.

{\em second part.}
${\cal K}^{f,{ I}}(x)\geq K(x)$.To prove this we first describe a certain class of {\em nice} balls covers.
Suppose that ${\alpha } =\{B_1(y_1,r_1),...,B_n(y_n,r_n)\}$ is a
ball cover of the metric space $X$  whose elements are balls with centers $y_i$
and radii $r_i$. We say that $\alpha$ is a {\em nice}
cover if $X\subset \stackunder {i}{\cup}B_i(y_i,\frac {r_i}{2} )$.
In other words  ${\alpha }$ is a {\em nice } cover if  dividing the radius of
the balls by $2$ we  have  again  a cover.

By Lesbegue Lemma about open covers and lemma \ref{16} we know that $\forall \lambda $ we can choose a rational $\epsilon $, such that there is
  a {\em nice } cover $\beta =\{B(I(s_1),2\epsilon),...,B(I(s_m),2\epsilon ) \}$    such that $| K^f
(x,\beta )-K^f (x)|<\lambda $   and thus we can consider covers of
the form $\beta $ to calculate Brudno's orbit complexity.
Now let us suppose to have a minimal length program  $p$ such that $ \forall i \in \{1,...,n\} ,d(U_i(p),T^i(x))< \epsilon$; there is an algorithm that calculates a sequence $\omega ^n$
relative to $\beta $. This algorithm can be coded in a program $p'$ with $|p'|=|p|+c$.  The algorithm calculates for each $i$ a string $z _i$
such that $I(z_i)=U_i(p)$ .
Then it finds a ball of $\beta ^{\prime }=\{B(I(s_1),\epsilon),...,B(I(s_m),\epsilon )\}$ that contains $I(z_i)$. This can be done
using the computable structure in the following way: for each $j,k$ it
calculates $D(z_i,s_j,k)$ until it finds a string $s\overline{_j}\in \{s_1,..,s_m\}$ such that
$d(I(z_i),I(s\overline{_j}))<\epsilon $ (this process must stop
because $\beta ^{\prime }$ is an open cover). The last statement
 implies that $
I(z_i)\in B(I(s_j),\epsilon )$ and $T^i(x)\in $ $B(I(s_j),2\epsilon) $.  Thus it is possible to construct a symbolic orbit relative to the cover $\beta $.$\Box$

Since ${ K}$ is independent of the choice of an interpretation we obtain

\begin{corollary}
If $X$ is compact, if $I$ and $J$ are computable interpretations (not necessarily from the same
computable structure) then ${\cal K}^{I,f}(x)={\cal K}^{J,f}(x)$.
\end{corollary}

If $X$ is not compact the two definition of orbit complexity are not equivalent. For example if $X={\bf R}$ with the standard computable structure, the trivially non chaotic dynamical system $T(x)=x+1$ is such that $K^{id}(x)=\infty,\forall x \in {\bf R}$ and ${\cal K}^{I,id}(x)=0\ \forall x \in {\bf R}$ equipped
with the standard computable structure.

If $X$ is not compact however the orbit complexity does not depends on the choice on $I\in {\cal I}$.

\begin{lemma}\label{20}
If $I ,J$ are computable interpretation functions from the same computable structure: $I ,J\in {\cal I}$ then ${\cal K}^{f,I}(x)={\cal K}^{f,J}(x)$.
So the orbit complexity does not depend on the choice of the interpretation $I$ in the computable structure $ {\cal I}$ and we can define ${\cal K}^{f,\cal I}(x)={\cal K}^{f,I}(x)$ for some $I \in {\cal I}$.
\end{lemma}
The proof is similar to the previous one and we omit it.

The  ``calligraphic'' orbit complexity is invariant for constructive isomorphisms of dynamical systems over non compact spaces, as it is stated in the following propositions.
We omit the proof that is similar to the previous ones.
\begin{theorem}
If $(X,d,T)$, $(Y,d',T')$ are topological dynamical systems over metric spaces
with  computable structures  ${\cal I},{\cal J}$  and 
$f$ is  onto and it is a  morphism between $(X,d,{\cal I})$ and $(Y,d',{\cal J})$
such that the following diagram 
\begin{equation}\label{equisopra}
\begin{array}{rcccl}
\ & \ & f & \ & \ \\
\ & X & \rightarrow & Y \ \\
T & \downarrow & \ & \downarrow & T'\\
\ & X & \rightarrow & Y \ \\
\ & \ & f & \ & \ \\
\end{array}
\end{equation}
\noindent commutes,if $x\in X$ and $y=f(x)\in Y$ then 
 ${\cal K}^{{\cal I},f}(x)\geq {\cal K}^{{\cal J},f}(y). $

\end{theorem}

\section{Initial condition sensitivity}

We now define two generalized indicators of initial condition sensitivity using 
a construction similar to the Brin-Katok local entropy.
The generalized indicators are constructed in a way that they have 
non trivial values in the weakly chaotic case. 
Let $X$ be a separable metric space  and $T$ a function $X\rightarrow X$.

Let us consider the following set:

$$ B(n,x,\epsilon)=\{ y \in X :d(T^i(y),T^i(x))\leq \epsilon \ \forall i \ s.t. \ \ 0\leq i \leq n\}.$$

$ B(n,x,\epsilon)$ is the set of points ``following'' the orbit of $x$ for $n$
steps at a distance less than $\epsilon$. As the nearby starting orbits of $(X,T)$ diverges  the set $ B(n,x,\epsilon)$ will be smaller and  smaller as $n$ increases. 
 The speed of decreasing of the size of this set considered
as a function of $n$ will be a measure of the sensitivity of the system
to changes on initial conditions.

Brin and Katok used the set $ B(n,x,\epsilon)$ for their definition 
of local entropy \cite{brinkatok}. In their paper the measure of the size of $ B(n,x,\epsilon)$ 
was the invariant measure of the set.

If we are interested to approximate the orbit of $x$ for $n$ steps we are interested to know how close we must approach the initial condition $x$ 
to ensure that the resulting approximate orbit is close to the orbit of $x$; another possible measure of the size of $ B(n,x,\epsilon)$ is then the radius of
the biggest ball with center $x$ contained in $ B(n,x,\epsilon)$.

$$ r(x,n,\epsilon)=\stackunder {B_r(x)\subset B(n,x,\epsilon)}{ sup\ r}.$$

Or the radius of the smaller ball that contains $ B(n,x,\epsilon)$ 

$$ R(x,n,\epsilon)=\stackunder {B_R(x)\supset B(n,x,\epsilon)}{ inf\ R}.$$

The generalized initial data sensitivity indicators at some point $x$ will be a family of real numbers indicating  the asymptotic behavior of how faster orbits coming  from a neighborhood of $x$ will diverge.
In a certain sense the behavior of $r(x,n,\epsilon)$ indicates the maximum
initial data sensitivity at $x$, while $R(x,n,\epsilon)$ indicates the
minimum initial data sensitivity at $x$.

 For this purpose we measure how faster
 $-log (r(x,n,\epsilon ))$ increases as $n$ increases, i.e. we consider the asymptotic behavior of the sequence  $-log (r(x,n,\epsilon ))$ as $n$ increases.
Let $f:{\bf N}\rightarrow {\bf R}$ be a monotone sequence, such that $\stackunder {n\rightarrow \infty}{lim}f(n)=\infty$. First we define 

\begin{definition}We define
$ r^f_{\epsilon}:X \rightarrow  {\bf R} $ as

$$ r^f_\epsilon (x)=\stackunder {n\rightarrow \infty}{limsup} \frac {-log( r(x,n,\epsilon))}{f(n)}$$

\noindent and

$ R^f_{\epsilon}:X \rightarrow  {\bf R} $ as

$$ R^f_\epsilon (x)=\stackunder {n\rightarrow \infty}{limsup} \frac {-log (R(x,n,\epsilon))}{f(n)}. $$

\end{definition}

The following lemma implies that $ r^f_\epsilon(x) $ and $ R^f_\epsilon(x) $ are  monotone functions with respect to $\epsilon$.

\begin{lemma}
If $\epsilon>\theta$ then $ r^f_{\epsilon}(x) \leq  r^f_{\theta}(x)$, $ R^f_{\epsilon}(x) \leq R^f_{\theta}(x)$.
\end{lemma}

Proof. Obvious $\Box$

\begin{definition}\label{d23}

We define the indicator of maximum initial data sensitivity at $x$ 
$  r^f(x):X\rightarrow {\bf R}$ as

$$ r^f(x)=\stackunder {\epsilon \in {\bf R}^+}{sup} r^f_{\epsilon}(x)$$

\noindent in the same way we define the indicator of minimum initial 
data sensitivity 
$ R^f(x):X\rightarrow {\bf R}$
$$ R^f(x)=\stackunder {\epsilon  \in {\bf R}^+ }{sup} R^f_{\epsilon}(x).$$

\end{definition}

The classical definition of dynamical system sensitive to initial 
conditions is related to our last  definitions. To say that a system is sensitive to initial conditions is equivalent to say that there is a $\delta$ such that  $r (x, n,\delta)$ is infinitesimal for all
the $x\in X$: the proof is immediate.

\begin{definition}
A dynamical system $(X,T)$ is said to have sensitive dependence on initial conditions
if there is a $\delta$ such that for each $x\in X$ and every neighborhood $U$ of $x$ there is $y\in U$ and $k\in {\bf N}$ such that $d(T^k(x),T^k(y))>\delta$. 
\end{definition}

 \begin{proposition}
A system has sensitive dependence on initial conditions if and only if
there is a uniform $\delta$ such that $\forall x \in X$ $r(x,n,\delta)$ goes to 0 as $n$ increases.
\end{proposition}

We give some example of different behaviors of $r(x)$ and $R(x)$ in
dynamical system over the interval $[0,1].$

\noindent {\em The identity map} $T(x)=x$. In this map $\forall n$ $B(n,x,\epsilon )=\{y\in
[0,1],|y-x|<\epsilon \}$ then if we choose for example $x=\frac 12$ we
have $R(\frac 12,n,\epsilon )=r(\frac 12,n,\epsilon )=\epsilon $.  Then $R^f_\epsilon(\frac 12)=r^f_\epsilon (\frac 12)=0$ for each $f$ (such that $f$ is not decreasing and $\stackunder {n \rightarrow \infty}{lim} f(n)=\infty $) and $R^f(\frac 12)=r^f(\frac 12)=0$.   The same arguments can be applied
to the  ergodic irrational translation on $[0,1]$: $T(x)=x+t \ (mod  \ 1 )$ where $t\notin {\bf Q}$ obtaining the same kind of initial data sensitivity as the identity (indeed both the maps are not sensitive to initial conditions).   

\noindent {\em The one dimensional baker's map} $T:[0,1]\rightarrow [0,1],T(x)=2x\ (mod\ 1).$ If we choose  for example $x=0$ , we
have  $B(n,0,\epsilon )=\{y\in [0,1], 0 \leq y \leq 2^{-n}\epsilon \},$ $R^f_\epsilon(0)=r^f_\epsilon (0)=\stackunder {n\rightarrow \infty}{limsup} \frac {n+log(\epsilon)}{f(n)} $ and $R^f(0)=r^f(0)=\stackunder {n\rightarrow \infty} {limsup} \frac {n}{f(n)} $ that is   $R^{id}(0)=r^{id}(0) =1 $, $R^f(0)=r^f(0)=0 $ if $f(n)=o(n)$ and $R^f(0)=r^f(0)=\infty $ if $n=o(f(n))$. In other words the baker's map has exponential initial data sensitivity
at 0.

\noindent {\em The piecewise linear Manneville map}  $T:[0,1]\rightarrow [0,1]$ defined in Equation \ref{**}.
 In this example (fig. 1) any neighborhood of the
origin $B_\epsilon =$ $[0,\epsilon )$ is subdivided in a sequence of
intervals $A_k=(\xi _k,\xi _{k-1}], \xi_k\sim k^{-\frac 1 {z-1}} $ and if $k>1$ then  $T(A_k)=A_{k-1}$. Let us
choose $x=0,$ then  $B(n,0,\xi _k)=[0,\xi _{k+n}).$ Then a point starting near to the origin goes away with a power law speed. By this we find $
r^f(0)=\stackunder {n\rightarrow \infty}{limsup} \frac {logn}{(z-1)f(n)}$ that is $r^f(0)=\frac 1{z-1}$ if $f(n)=log(n)$.  In other words the map $T(x)$  has power law
sensitivity to initial condition at the origin. 
In \cite{ionew} it is proved that (while the sensitivity to initial condition at the origin is a power law)  for almost all other  points in $[0,1]$ we
have a stretched exponential sensitivity.

Let us consider the following map $ T:[-1,1]\times [-1,1]\rightarrow [-1,1]\times [-1,1] $  defined as 
$T\left( \begin{array}{c}x \\ y \end{array} \right)=\left( \begin{array}{c}
x +y  \ (mod \ 2)   \\ 
y \ (mod \ 2)  \end{array}
\right) . $ 
This example is a simplification of the Casati map \cite{casati}.
Let us consider the initial data sensitivity at the origin.
It can be  calculated that $r(0,n,\epsilon )$ decreases asymptotically (modulo
lower order terms)
 as $\frac {Const} n$  then $r^{log (n)}(0)=1$ and $r^f(0)=\infty$ for any $o(f(n))=log(n)$.
Conversely $R(0,n,\epsilon)$ is constant because the map 
is the identity on the axis $y=0$  and then $\{y=0\}\cap B(0,\epsilon)\subset B(0,n,\epsilon)$ for each $n$.
From this, since $R(0,n,\epsilon) $ is not infinitesimal it follows that $R^f(0)=0$ for each $f$.

\section{Complexity of points}
Now we define a function ${ S}^I(x,\epsilon ):X\times {\bf R}\rightarrow {\bf R}$, the function is a measure of the complexity of the points of $X$. The function is not increasing with respect to $\epsilon $  and  measures how much information is necessary to approximate a given point of $X$ with accuracy $\epsilon$. Thus it is a function that does not depend on the dynamics.
  In \cite{io1} a definition of local entropy for points of metric spaces was based on this idea and connections between $S$ and the concept of dimension are shown. In particular $S(x,\epsilon)$ is related to the local dimension of $X$ at $x$. 

We now recall the results and reformulate it in a way to be used in dynamical systems.

\begin{definition}If $(X,d,{\cal I}) $ is a metric space with a computable structure $\cal I$,  $ I$ is an interpretation function in $\cal I$, ${\cal U}$ an universal computer. The information contained in the point $x$ with respect to
the accuracy $\epsilon$ is defined as: 
\begin{equation} \label{1} 
{ S}^I(X,x,\epsilon )={\min} \left\{ |p|\ \ s.t.\ \ d(I({\cal U}(p)),x)<\epsilon\right\}.
\end{equation}
\end{definition}

In the notation ${ S}^I(X,x,\epsilon )$, in the following the space $X$ will
not be specified when it is clear from the context. 

The function $S$ depends on the interpretation $I$,
the function $S$ depends also on the choice of ${\cal U}$.
Changing the choice of the universal Turing machine $\cal U$ will change the value of $S$ by a constant.
Since in the following we are interested to the asymptotic behavior of 
$S^I(X,x,\epsilon )$ when $\epsilon $ goes to $0$ the choice of $\cal U$
is not relevant.

%
%

Now we see the relations between $S$ and the concept of  dimension.
Let $F\subset X$ be a completely bounded metric subspace of $X$.
We recall the definition of box counting dimension:
let 
\begin{equation}\label{n} n_{\delta}(F)=
 min\{ n\in {\bf N},such \ that \ \exists x_1,...,x_n \ such \ that \stackunder i\cup B(x_i,\delta)\supset F \} 
\end{equation}
 be the minimum number of $\delta$ balls that cover $F$, the upper and
lower box counting dimension of $F$ are defined as:
$$
\underline {d}_B (F)=\stackunder {\delta \rightarrow 0}{liminf} \frac {log (n_\delta (F))}{-log \delta}
$$

$$\overline {d}_B (F)=\stackunder {\delta \rightarrow 0}{limsup} \frac {log (n_\delta (F))}{-log \delta}.$$

An equivalent definition can be given by considering instead of $n_\epsilon$ the maximum number of disjoint balls that are contained in $F$.
If $\{x_0,...,x_i\}\subset F$ are such that $d(x_i,x_j)>\epsilon $ $ \forall i,j,i\neq j $ we say that the $\{x_0,...,x_i\}$ forms an $\epsilon$ net.

Let us consider
the maximum cardinality ${n'}_\delta (F)$ of a $\delta$ net with centers in $F$.
It can be easily proved (see e.g. \cite{falconer} pag 39) that $n_{4\delta} (F)\leq {n'}_\delta (F)$ and  $ {n'}_\delta (F) \leq  n_\delta (F)$ and then the two quantities give rise to the same definition of box counting dimension.

The box counting dimension is widely used in dynamical systems and their
applications. Some of its feature differs from the Hausdorff dimension.
For example the following holds 
\begin{proposition}\label{31} If $X\subset F$ then
$$\overline {d}_B(X)=\overline {d}_B(\overline X),\ \underline {d}_B(X)=\underline {d}_B(\overline X)$$ where $\overline X$ denotes the closure of $X$.
\end{proposition}

\noindent for the proof see for example  \cite{falconer} page 44. Hence countable sets may have positive dimension. In general we have $d(X)\leq \underline{d}_B(X)\leq \overline{d}_B(X) $, where $d(X)$ is the Hausdorff dimension of $X$.


\begin{proposition}\label{32}
Let $(Y,d,{\cal I})$ is a metric space with a computable structure $I\in {\cal I}$   and $ \overline{d}_B(Y)=d $ 
then  for each $ x$ and $\delta $,  if $\epsilon$ is
small enough, i.e. $\exists \overline{\epsilon}$ such that $\forall \epsilon <\overline{\epsilon}$ $$ S^I(x,\epsilon)\leq -(d+\delta)log(\epsilon).$$
\end{proposition}
{\em Proof.} 

We prove that $\exists C$ such that for each $x \in X$ and $\delta >0$ exists $\overline {\epsilon} $ such that $\forall \epsilon< \overline{\epsilon}$ $$ S^I(x,\epsilon)\leq -(d+\frac\delta 2)log(\epsilon)+log(log(\epsilon))+C.$$ from this proposition \ref{32} follows easily.

First let us suppose that $\epsilon $ is of the form $\epsilon = 2^{-z},z\in {\bf N}$.
We describe a recursive procedure $P(n_x,\epsilon ):{\bf N}\times {\bf Q} \rightarrow \Sigma $
such that each point $x$ of $Y$ is  $3\epsilon$ near to 
$I(P(n_x,\epsilon ))$ for some $n_x \in {\bf N} $ (i.e. $\forall x \in X \ \exists n_x  \ s.t. \ d(x,I(P(n_x,\epsilon )))<3\epsilon$).
 $P(n_x,\epsilon)$ builds an  $\epsilon$ net  of $n_x$ balls and stops
giving the center of the last ball found as the output.
Let $s_0$ be the first string of $\Sigma $ lexicographic ordered. 

 $P(n_x,\epsilon)$ starts  from $I(s_0)$ with a list $L$ containing $s_0$.
Then it consider each string of $\Sigma$ (in lexicographic order) $s_0,s_1,...$
at the first step $L=\{s_0\}$, the procedure looks for a string $s_{i_1}$
such that $$d(I(s_{i_1}),I(s_0))\geq \epsilon$$ (as it was explained in section 2.1) then it adds $s_1 $ to the list.
And so on, by induction at the $n-$th  step it searches for  $s_{i_n}$ such that $$\forall s \in L \ d(I(s_{i_n}),I(s))\geq \epsilon
$$ and it add it to $L$. When $n=n_x $ the procedure stops and outputs
the last string in the list. 

In this way we build an $\epsilon$ net.
Let us call $n_P$ the maximum number of centers that can be found by such a 
procedure, i.e. if $n_x\leq n_P$ then the procedure stops.

Since the dimension of $Y$  is $d$ and at most we can have ${n'}_\epsilon$ balls. For each $\delta >0$, $\exists \overline{\epsilon}$ such that $\forall \epsilon <\overline{\epsilon}$ $$n_x\leq n_P \leq {n'}_\epsilon \leq  \epsilon ^ {-(\overline{d}+\frac \delta 2)}$$

Now we remark that each  $x\in Y$ must be $3\epsilon$ near to one of these centers $x_0,...,x_{n_P}$, because if there is $\tilde {x}$, such that$\forall i \leq n_P$ $d(\tilde x,x_i)>3\epsilon$ then the procedure will find it and add it to the net, contradicting the maximality of $n_P$.
Then the information that is sufficient to approximate $x$ at accuracy $\epsilon$ is: the code for the procedure $P$ (which is constant with respect to $\epsilon$), the number $n_x\leq \epsilon ^{-(d+\frac \delta 2)}$, whose binary  length is
$\leq log (n_x) +1$, and the information that is necessary to give $\epsilon =2^{-z}$ that is about $log(z)=log(log(\epsilon ))$.

Then for each $ \delta$ we have for $\epsilon $ small enough of the form $\epsilon =2^{-z}$
\begin{equation}\label{fascino1}
 S(x,\epsilon)\leq log(n_x)+log(log(\epsilon))+C\leq log ({(\frac{\epsilon}3})^{-(d+\frac \delta 2)})+log(log(\epsilon ))+C \leq 
\end{equation}
$$ \leq  {-(d+\frac \delta 2)}log (\epsilon) +log(log(\epsilon ))+C.
$$
 
If $\epsilon $ is not of the above form the statement follows from the remark that if $\epsilon $ is small enough and \ref{fascino1} holds  even for $ 2^{[log \epsilon]}$  then
$$ S^I(x,\epsilon )\leq  S^I(x,2^{[log \epsilon]} )\leq (d+\frac \delta 2)log(2^{[log \epsilon]})+log(log(\epsilon ))+C\leq $$
 
$$\leq  (d+\frac \delta 2)(log(\epsilon)+1)+log(log(\epsilon ))+C. $$
$\Box$

\noindent {\bf Remark. }If $T$ is continuous then  the closure of an orbit is an invariant set.

\begin{lemma}\label{33}
If $x_n$ and $y_n$ are sequences in $X$ and $d(x_n,y_n)\leq 2^{-n}$  then 
$\overline {d}_B(\{ x_i\}_{i\in{\bf N}} )=\overline {d}_B(\{ y_i\}_{i\in{\bf N}} )=\overline {d}_B(
\overline {\{ x_i\}_{i\in {\bf N}}\cup\{ y_i\}_{i\in{\bf N}}}).$ 
\end{lemma}
{\em Proof.} 
Let $\{B_i\}$ be a minimal cover of $\{ x_i\}$ made of balls $B_i$ with radius $\epsilon $ and $\#\{ B_i \}=n_\epsilon $ as defined in \ref{n}, then $\{ 2B_i\}$\footnote{$2B_i$ is a ball with the same center as $B_i$ and double radius.} is a cover of all points of $\{y_i\}$ but some of the $y_i$ such that $d(x_i,y_i)\geq \epsilon$ and these points 
are at most $- log \epsilon $ many. Then 
$ n_{2\epsilon}(\{y_i \}) \leq n_{ \epsilon } (\{ x_i  \}) - log \epsilon $,
then $$ \overline {d}_B(\{y_i\})=\stackunder{\epsilon \rightarrow 0}{limsup} \frac {log(n_{2\epsilon}(\{ y_i\}))}{log 2\epsilon}\leq \stackunder {\epsilon \rightarrow 0}{limsup} \frac {log(n_{\epsilon}(\{x_i\})-\log \epsilon )}{log 2\epsilon}=\overline {d}_B(\{x_i\})$$
 and then $\overline {d}_B(\{ y_i\})\leq \overline {d}_B(\{x_i\})$.
In the same way it is possible to obtain the reverse inequality and by the use of  
proposition \ref{31}
the proof of $\overline {d}_B(\{ x_i\} )=\overline {d}_B(\overline {\{ x_i\}\cup\{ y_i\}} )$. $\Box$

\begin{proposition}\label{34}
 Let $(X,{\cal I})$ be a metric space with a computable structure and $T:X\rightarrow X$  such that $(X,T) $ is constructive. If $x$ is contained in the closure of the orbit of an ideal point: $ x\in O=\overline {orb(I({s_0}))}$ then $\exists I\in {\cal I}$, such that for each $\delta $ if $\epsilon $  is small enough
 $$S^I(x,\epsilon )\leq -(\overline {d}_{B}(O)+\delta)log \epsilon.$$
 \end{proposition}
{\em Proof.}
Let us consider the following interpretation function on $X$
$I(s)= \left\{ \begin{array}{c}
I_{0}(s') \ if s=0s'\\
I_{1}(s') \ if s=1s'
\end{array}\right.    $
where $ I_0$ is some interpretation in ${\cal I}$  and $I_1 $ is constructed
 in the following way: if $A(n,\epsilon,s_0)$ is the algorithm to follow the orbit of $I(s_0)$ introduced in Lemma \ref{lemma12}, then 
$$I_1(s)=I_0(A(s_0,c(s),c(s)))$$
where $c(s)$ is the natural number  associated to   $s$ (Section. 2). If $T$ is constructive then $I$ and $I_0$ are equivalent. Since the distance of the approximate orbit $I_0(A(s_0,c(s),c(s)))$  from the real orbit of $x$ 
decreases exponentially,
 by the above lemma \ref{33} $$\overline {d}_{B}(\overline { \{I_1(\Sigma)\}})=\overline {d}_{B}(O)$$
 and if we consider the space $\tilde {O} =\overline {\{I_1(\Sigma)\}}\cup O$ we have ${\overline d}_{B}(\tilde {O})={\overline{d}}_{B}({O})$. Now  we consider the metric space $\tilde {O}$
 with the interpretation $I_1$ we have trivially that $\exists C$ such that 
$$ S^I(X,x,\epsilon)\leq S^{I_1}(\tilde {O},x,\epsilon )+C $$ where $S(X...)$ means the information need to approximate $x$ considered as a point of $X$. Moreover by Proposition \ref{32} for each $\delta$ if $\epsilon$ is small enough $ S^{I_1}(\tilde {O},x,\epsilon) \leq -({\overline {d}}_{B}(\tilde {O})+\delta)log (\epsilon )  $,  and then $S^I(X,x,\epsilon)\leq S^{I_1}(\tilde {O},x,\epsilon )+C \leq  -({\overline {d}}_{B}(\tilde {O})+\delta)log \epsilon  +C=-(\overline d _B (O)+\delta)log (\epsilon )+C $ 
and  the statement
is proved.$\Box$

\begin{proposition}\label{qeu}
 For each $I$, for each $x$ but a set $F$  with $d$ dimensional Hausdorff measure $H_d(F)=0$, for each $\delta$, if $  \epsilon $ is small enough
$$S^I(x,\epsilon)\geq -(d - \delta )log \epsilon . $$
\end{proposition}
{\em Proof.} 
By theorem 12 \cite {io1}, page 1293 the set \begin{equation}\label{mistero}
W^d=\{ x\in X\ s.t. \stackunder {i\in {\bf N}, i\rightarrow \infty }{liminf}\frac{S^{ I}(x,2^{-i})}{i}\leq d\}
\end{equation} has dimension $\leq d$.

By contradiction if $x$ is such that there is a $\delta$ such that there is a sequence $\epsilon _i \rightarrow 0 $ such that for each $i$  $S(x,\epsilon _i)< -(d - \delta )log( \epsilon _i)$,
then $\frac{S(x,\epsilon _i)}{-log(\epsilon _i)}<d-\delta $,
then $$\stackunder {\epsilon \rightarrow 0}{liminf}  \frac{S(x,\epsilon )}{-log\epsilon}\leq d-\delta$$ if $\epsilon=2^{-\gamma}$ then $\stackunder{\gamma \rightarrow \infty} {liminf} \frac{S(x,2^{-\gamma} )}{\gamma}\leq d-\delta $. Since 
$S(x,2^{-\gamma-1}) \leq S(x,2^{int(-\gamma)})\leq S(x,2^{-\gamma+1})$
then $\stackunder {i\in {\bf N},i \rightarrow \infty}{liminf}  \frac {S(x,2^{-i} )}{i}\leq d-\delta $ and then by \ref{mistero} $x \in \stackunder {d'<d}{\cup } W^{d'}$
 and this union of sets has $d$ dimensional Hausdorff measure $0$  $\Box$

{If $\mu $ is a measure on  $X$ let us consider }

{$$\underline { d}_{\mu}(x)=\stackunder {\epsilon \rightarrow \infty}{{liminf}}\frac {log(\mu(B(x,\epsilon)))}{log(\epsilon )}$$}
 ${ d}_{\mu}(x)$ is a sort of local dimension of the measure $\mu$ with respect
to the metric of $X$ at the point $x$ and was widely studied in the 
theory of dynamical systems (see for example the book \cite{pes}).

\begin{proposition}\label{lemmone} For $\mu $ almost each $x$, for each $\delta $, $\exists \overline {\epsilon } $ s.t. $\epsilon\leq \overline {\epsilon } $ implies 
$$S(x,\epsilon)\geq -(\underline {d}_{\mu}(x)-\delta ) log \epsilon . $$
\end{proposition}

For the proof of this proposition we need some additional lemmata and remarks.
\begin{remark}\label{magia} By the proof of Proposition \ref{qeu} 
if we have that  $x$ is such that $\exists \delta$ and $\epsilon _i\rightarrow 0$  $s.t.$ $S(x,\epsilon _i)< -(\underline {d}_{\mu}(x)-\delta ) log (\epsilon _i)$, 
then $$S^{exp}(x):=\stackunder {i\in  {\bf N},i\rightarrow \infty }{liminf} \frac {S(x,2^{-i})}i \leq \underline {d}_\mu (x)-\delta.$$ 
\end{remark}

\begin{lemma}
Let $d,c,\delta >0$, let $B\subset X$ such that $\mu (B)>0$
and $\forall x\in B \ \underline {d}_{\mu} (x) 
\geq  d+c  $. If $A=\{x\in B \ s.t. \ \exists \delta  \ s.t. \ S^{exp}(x) \leq d-\delta \}$, then $\mu(A)=0$. 
\end{lemma}

{\em Proof. }Since $\forall x\in  A $ $\underline {d}_{\mu}(x) \geq d+c $ then $\forall x\in A $,  if $\gamma $ is big enough,  $$\mu(B(x,2^{-\gamma}))<2^{-\gamma(d+\frac c2 )}$$
 then there is a $\overline {\gamma}>0$ and a set $A'\subset A$ with $\mu (A')>0$ such that if $x\in A'$ $\forall \gamma >\overline {\gamma} \ \mu(B(x,2^{-\gamma}))<2^{-\gamma (d+\frac c 2)}$ uniformly on all $A'$.

By the remark above for each $k$, each $x\in A$ is contained in some ball
$B(I({\cal U}(p)),2^{-i})$ with $|p|>k$ and $\frac {|p|}i<d$, that is 
 $\forall k\in {\bf N}$ $A\subset \stackunder { \{ p|k \leq |p|\} } {\cup} B(I({ \cal U}(p)),2^{-\frac {|p|}{d}})$.
Let us consider the measure of one of these balls . If $\frac {|p|}{d} -1>\overline \gamma  $ and $B(I(s),2^{- \frac {|p|}{d}})$, is such that $B(I(s),2^{- \frac {|p|}{d}}) \cap A\neq \emptyset$ then $B(I(s),2^{- \frac {|p|}{d}})\subset B(y\in A,2^{\frac {-|p|}d+1})$ then exists $y\in A$ such that  $\mu(B(I(s),2^{- \frac {|p|}{d}}))\leq \mu( B(y,2^{\frac {|p|}d+1}))\leq 2^{-\frac {|p|+d}d (d+\frac c2)}.$
Then, since for each $k$ there are at most $2^{k+1}$ programs of length $k$
  $$\mu (A)\leq \stackunder {\{ p|k\leq |p|,B(I(s),2^{- \frac {|p|}{d}}) \cap A\neq \emptyset \}}{\sum}\mu ( B(I({\cal U}(p)),2^{-\frac {|p|}d})) \leq \stackunder {i\geq k}{\sum} 2^{i+1}2^{-\frac {d+\frac c2}di -d -\frac c2 } $$ since $i -\frac {d+\frac c2}di <0 $ and $k$ is arbitrary  the last sum can be set as small as wanted the statement is proved.  $\Box$

The last lemma implies an a.e. relation between dimension and $S(x,\epsilon)$:

\begin{lemma}
For $\mu $ a.e. $x$ 
$S^{exp}(x) \geq \underline{d}_\mu(x)$.

\end{lemma}

{\em  Proof.} If   conversely  $ S^{exp}(x)< \underline{d}_\mu(x)$
on a set $A'$ with $\mu(A')>0$ it is possible to find a constant $c$ and a set $A''$, $\mu(A'')>0$ such that $S^{exp}(x)<c< \underline{d}_\mu(x)$ on $A''$, by lemma  16 we obtain  $\mu(A'')=0$. $\Box$

By this last lemma  and Remark \ref{magia} the proof of Proposition \ref{lemmone} follows directly.

\section{Initial data sensitivity and orbit complexity} 

Now we are ready to state the first proposition linking orbit complexity
to initial data sensitivity, for the proof see \cite{ionew}.

\begin{proposition}\label{quellagrossa}
 If $(X,T)$ is a dynamical system on a space $(X,d,{\cal I})$ with a computable structure $\cal I$, $I\in {\cal I}$ and $T$ is constructive. There are constants $c_1 $ and $ c_2$ such that for all $x\in X,n\in {\bf N},\epsilon \in {\bf R}^+$ 
\begin{equation}\label{quella1}
{\cal K}^I(x,n,2\epsilon)<{S}^I(x,r(x,n,\epsilon))+log (n)+c_1
\end{equation}

\begin{equation}\label{quella2}
{ S}^I(x,R(x,n,3\epsilon))\leq {\cal K}^I (n,x,\epsilon)+c_2.
\end{equation}
\end{proposition}

We remark that the presence of the term  $''+log(n)''$ in equation \ref{quella1} is due to our simple definition of AIC. If instead of the plain algorithmic
information content of a string $s$ we had used $AIC(s|lenght(s))$ (the 
algorithmic information content of a string given its length) the term
$''+log(n)''$ would disappear from Equation \ref{quella1}.

In the next theorems we prove some relations between our generalized indicators of initial 
data sensitivity, orbit complexity and various notions of  dimension of the space $X$.

In next theorem it is  also considered the dimension of the closure of an orbit
containing our initial point $x$. This gives a relation that in a certain  sense involves
the dimension of attractors.

\begin{theorem}
If $(X,T)$ is a constructive dynamical system over  a compact metric space with a computable structure $(X,d,{\cal I})$. If  $f$ is such that $log(n)=o(f(n))$
then 
\begin{equation}\label{brivido} {K^f(x)}\leq  \overline{d}_B(X){ r^f(x)}.
\end{equation} 

If $f(n)=log(n)$ then we have $$ K^{log}(x) \leq  \overline{d}_B(X){ r^{log}(x)}+1.$$

If $x$ is contained in the closure of the orbit $ \overline {O_{I(s_0)}}$ of some ideal point $I(s_0)$ and  $log(n)=o(f(n))$ then
 $$ {K^f(x)}\leq \overline{d}_B(\overline {O_{I(s_0)}}){ r^f(x)} $$
 and  $$ K^{log}(x)\leq \overline{d}_B(\overline {O_{I(s_0)}}){ r^{log}(x)}+1.$$ \end{theorem}

{\em Proof.} Since the space is compact then ${\cal K}^{f,I}(x)=K^f(x)$ and  then it is sufficient to prove the statement for  ${\cal K}^{f,I}(x)$. Let us  suppose that $r(x,n,\epsilon)$ is infinitesimal. For each $\delta >0$ if $n$  is big enough, ( and then $r(x,n,\epsilon)$ is  eventually small),
by proposition \ref{quellagrossa} and proposition \ref{32}

$$\frac {{ \cal K}^I(x,n,2\epsilon)}{f(n)}<\frac{{S}^I(x,r(x,n,\epsilon))+log (n)+c_1}{f(n)}\leq $$$$ \leq \frac{-(\overline{d}_B(X)+\delta) log(r(x,n,\epsilon)) +log (n) +C}{f(n)}.$$

If $log(n)=o(f(n))$ by taking the limsup for $n\rightarrow \infty$ and limit for $\epsilon \rightarrow 0$

$$ {{\cal K}^f(x)}\leq (\overline{d}_B+\delta){ r^f(x)}$$

since $\delta $ is arbitrary we obtain 
${{\cal K}^{f,{\cal I}}(x)}\leq  \overline{d}_B(X){ r^f(x)}$.

If $f(n)=log(n)$ then by Proposition \ref{quellagrossa} we obtain in the same way as before $$\frac {{ \cal K}^I(x,n,2\epsilon)}{log(n)}<\frac{{S}^I(x,r(x,n,\epsilon))+log (n)+c_1}{log(n)}$$ and we conclude as before.

Now suppose that $x$ is contained in the closure of the orbit of some ideal point. By Proposition \ref{34} we have that there is an $ I\in {\cal I}$ and $ C\in {\bf N}$, such that for each $\delta $ if $\epsilon $ is small enough 
 $$S^I(x,\epsilon )\leq (\underline {d}_{B}(\overline {O_{I(s_0)}})+\delta)log \epsilon +C.$$ 
Now as before by Proposition \ref{quellagrossa} 
$$\frac {{\cal K}^I(x,n,2\epsilon)}{f(n)}<\frac{{S}^I(x,r(x,n,\epsilon))+log (n)+c_1}{f(n)}\leq $$ $$  \leq \frac{-(\overline{d}_B+\delta) log(r(x,n,\epsilon)) +log (n) +C}{f(n)}$$
for each $\delta >0$ if $n$ is big enough  and we can conclude as before.
If $r(x,n,\epsilon)$ is not infinitesimal (and $\forall f$ $r^f(x)=0$) then by  Proposition \ref{quellagrossa} we have that ${\cal K}^{I}(x,n,\epsilon)\leq \log n +C$ then the statement follows 
trivially. 
 $\Box$ 

We remark that in the  above theorem the difference between
the statements when $f(n)=log(n)$ and
$log (n)=o(f(n))$ comes from the definition of the information content
of a string (see remark before proposition \ref{quellagrossa}).

%
%
\begin{theorem} If $X$ is compact, $(X,T)$ is constructive, 
for each $x$ but a set $F$ with  $d$ dimensional Hausdorff measure $H_d(F)=0$
$$ {K^f(x)}\geq d\ {R^f(x)}.
$$
\end{theorem}
{\em Proof.}
Again the space is compact and then ${\cal K}^{f,I}(x)={\cal K}^f(x)$ then it is sufficient to prove the statement for  ${\cal K}^{f,I}(x).$
By proposition \ref{quellagrossa} 
$$
 { S}^I(x,R(x,n,3\epsilon))\leq {\cal K}^I (n,x,\epsilon)+c_2.
$$ If $R(x,n,\epsilon )$ goes to $0$ as $n$ goes to infinity
by Proposition \ref{qeu}
 we obtain that there is a set $F$ with $H_d(F)$ such that if $x\in X-F$ $\forall \delta$ , if $n$ is big enough

$$  -(d-\delta ) log (R(x,n,\epsilon))\leq  {{ S}^I(x,R(x,n,3\epsilon))}\leq {{ \cal K}^I (n,x,\epsilon)+c_2}.$$

Dividing by $f(n)$ and  taking the limsup for $n\rightarrow \infty$ and limit for $\epsilon \rightarrow 0$ as before we have
$$ {{\cal K}^f(x)}\geq (d-\delta ) {R^f(x)}
$$ and we conclude as before.
If  $R(x,n,\epsilon )$ is not infinitesimal then $R^f(x)=0$ and the statement
is trivially true. $\Box$

\begin{theorem}
If $\mu $ is a measure on $X$ then for almost each $x$

$$  {K^f(x)}\geq  \underline {d}_\mu(x)  {R^f(x)}
$$
\end{theorem}

{\em Proof.} As before
by Proposition \ref{quellagrossa} and Proposition \ref{lemmone}, for almost each $x\in X$
$$
 -(\underline{d}_{\mu(x)}-\delta )log (R(x,n,\epsilon))\leq { S}^I(x,R(x,n,3\epsilon))\leq {\cal K}^I (n,x,\epsilon)+c_2
$$ and the proof can be concluded as above. $\Box $

\begin{remark}We remark that if $X$ is  totally bounded but
not compact all the above theorems hold with ${\cal K}^{f,I}(x)$ instead of $K^f(x)$. 
\end{remark}

We finally remark a connection with the theory exposed above and 
some applicative purposes. The algorithmic information content of a string
is not a computable function, it is impossible to calculate $AIC(s)$ for each $s$ by an algorithm.
By this 
 then the information content of an orbit as it is defined in this
paper 
is not computable. 

It follows that the concept of orbit complexity as it is defined in this paper
is not a quantity that can be experimentally measured.

To solve this problem
let us consider a data compression algorithm.
Since in the compressed string there is all the information that is 
necessary to reconstruct the original string
an 'approximate' measure for the information content of the string is the length of the string after it is compressed 
by this coding procedure.
This approximate measure is as accurate as efficient is the data compression procedure.

Using compression algorithms instead of the algorithmic information content it is possible to obtain a 'computable' notion of orbit complexity. It can be proved that  in the positive entropy case the computable orbit complexity is a.e. equivalent to the AIC based
one (\cite{io3},\cite{Licata}). 
Such a definition of computable orbit complexity allows numerical
investigations about the complexity of unknown systems. Unknown systems underlying  for example some given time series or experimental data.

In \cite{ZIPPO} and \cite{menconi} a particular compression algorithm
is developed to apply these ideas to the weakly chaotic case and 
numerical investigations are performed by directly measuring the information content of the orbits of some weakly chaotic 
system (e.g. the Manneville map and the logistic map at the chaos threshold
). The results agree with the theory exposed in this paper.

In our opinion the existence of a computable version of the orbit complexity
motivates from the applicative point of view the study of the  orbit complexity itself and its relations between the other  measures of  the chaotic behavior of a system.

\end{document}